\nonstopmode
\documentclass[a4paper, 10pt, reqno]{amsart}
\usepackage{latexsym}
\usepackage{fancyhdr}
\usepackage{amsmath, amssymb}
\usepackage[ansinew]{inputenc}
\usepackage[all]{xy}
\usepackage{verbatim}
\theoremstyle{plain}
\newtheorem*{lemma*}{Lemma}
\newtheorem{lemma}[subsection]{Lemma}
\newtheorem*{theorem*}{Theorem}
\newtheorem{theorem}[subsection]{Theorem}
\newtheorem*{proposition*}{Proposition}

\newtheorem*{corollary*}{Corollary}
\newtheorem{corollary}[subsection]{Corollary}
\theoremstyle{definition}
\newtheorem*{definition*}{Definition}

\newtheorem*{example*}{Example}
\newtheorem{example}[subsection]{Example}
\newtheorem*{remark*}{Remark}
\newtheorem*{remarks*}{Remarks}
\newtheorem{remark}[subsection]{Remark}

\numberwithin{equation}{section}
\def\x{\times}
\def\ox{\otimes}
\def\al{\alpha}
\def\be{\beta}
\def\ga{\gamma}
\def\de{\delta}

\def\la{\lambda}

\def\si{\sigma}

\def\vh{\varphi}

\def\om{\omega}
\def\Ga{\Gamma}
\def\De{\Delta}

\def\La{\Lambda}

\def\Ph{\Phi}

\def\Om{\Omega}

\def\R{\mathbb{R}}

\def\p{\partial}

\def\<{\langle}
\def\>{\rangle}
\def\-{\backslash}
\renewcommand{\o}{\circ}

\let\on=\operatorname
\def\sr#1%
{\ifmmode{}^\dagger\else${}^\dagger$\fi\ifvmode
\vbox to 0pt{\vss
 \hbox to 0pt{\hskip\hsize\hskip1em
 \vbox{\hsize3cm\eightpoint\raggedright\pretolerance10000
 \noindent #1\hfill}\hss}\vss}\else
 \vadjust{\vbox to0pt{\vss%
 \hbox to 0pt{\hskip\hsize\hskip1em%
 \vbox{\hsize3cm\eightpoint\raggedright\pretolerance10000%
 \noindent #1\hfill}\hss}\vss}}\fi%
}

\title[On the continuous cohomology of diffeomorphism groups]
{On the continuous cohomology of diffeomorphism groups}

\author[M.V.\ Losik ]
{M.V\ Losik}
\address{M.V.~Losik: Saratov State University, 
Astrakhanskaya, 83, 410026 Saratov, Russia}
\email{losikmv@info.sgu.ru}

\begin{document}
\begin{abstract} Suppose that $M$ is a connected orientable $n$-dimensional manifold and $m>2n$. If  $H^i(M,\R)=0$ for $i>0$, it is proved that for each $m$ there is a monomorphism $H^m(W_n,\on{O}(n))\to H^m_{\on{cont}}(\on{Diff}M,\R)$. If $M$ is closed and oriented, it is proved that for each $m$  there is a monomorphism 
$H^m(W_n,\on{O}(n))\to H^{m-n}_{\on{cont}}(\on{Diff}_+M,\R)$, where $\on{Diff}_+M$ is a group of preserving orientation diffeomorphisms of $M$.    
\end{abstract}
\date{\today} 
\maketitle
\section{Introduction}
Let $M$ be a connected orientable $n$-dimensional manifold. 
Denote by $D$ the group $\on{Diff}M$ of diffeomorphisms of $M$ and by $D_+$ the group of diffeomorphisms of $M$ preserving the orientation of $M$ whenever $M$ is oriented.  Later we consider $D$ or $D_+$ as an infinite-dimensional Lie group by \cite{KM} and use the calculus of differential forms and vector fields on infinite-dimensional manifolds developed in this book. 

In the present paper we study the continuous cohomology $H^*_{\on{cont}}(D,\R)$ and $H^*_{\on{cont}}(D_+,\R)$ of the groups $D$ and $D_+$ with values in the trivial $D$ and $D_+$-module $\R$. The main known result about the cohomology $H^*_{\on{cont}}(G_+,\R)$ is Bott's theorem (\cite{B}) about $(n+1)$-cocycles 
on the group $D_+$ for the closed oriented $M$. These cocycles are obtained from some $(2n+1)$-cocycles of the complex $C^*(W_n,\on{O}(n))$ of relative with respect to the group $\on{O}(n)$ cochains of the Lie algebra of formal vector fields $W_n$  and are expressed via integrals using a Riemannian metric on $M$. Moreover, for $M=S^1$ it is known (\cite{F}) that 
the cohomology $H^*_{\on{cont}}(D_+,\R)$ is a ring with two generators 
$\al,\be\in H^2(D_+,\R)$ which satisfy the only relation $\be^2=0$. 
 
Next, for brevity, put $H^p(D)=H^p_{\on{cont}}((D,\R)$, $H^p(D_+)=H^p_{\on{cont}}(D_+,\R)$, and $H^p(M)=H^p(M,\R)$. Let $H^i(M)=0$ for $i>0$. We prove that in this case, for each $m>2n$, there is a monomorphism 
$H^m(W_n,\on{O}(n))\to H^m(D)$. Let $M$ be closed and oriented. We prove that in this case, for each $m>2n$, there is a monomorphism 
$H^m(W_n,\on{O}(n))\to H^{m-n}(D_+)$. In particular, in the last case Bott's 
cocycles define a part of the monomorphism above for $m=2n+1$. 

The main idea of the proof is the following. Denote by $\Om^*(M)$ the de Rham complex of the mahifold $M$ and by $\on{Vect}M$ the Lie algebra of vector fields on $M$. First we interpret $H^*(D)$ in terms of some double complex. Then we  compare the cohomology of this double complex with the diagonal cohomology of the double complex $C^*_\De(\on{Vect}M;\Om^*(M))$ of the  Lie algebra $\on{Vect}M$ with values in the $\on{Vect}M$-module $\Om^*(M)$ and use the known facts on this cohomology. 

Section \ref{prel} contains the main algebraic constructions which will be used afterwards. Namely, one introduces a double complex $C^*_{\on{cont}}(G;K)$ for a topological group $G$ acting on a topological cochain complex $K$, one  recalls  some results on the cohomology $H^*(W_n)$ of the Lie algebra $W_n$ of formal vector fields and the relative cohomology $H^*(W_n,\on{O}(n))$,  the space $S(M)$ of frames of infinite order of $M$ and the canonical Gelfand-Kazhdan form with values in $W_n$ on $S(M)$. 

In section \ref{Dcocycles} one proves the theorem showing how to construct the cocycles of the groups $D$ and $D_+$ from cocycles of the complex $H^*(W_n,\on{O}(n))$ (Theorem \ref{group cocycle}, corollaries \ref{cocycle1} and \ref{cocycle2}). Moreover, one recalls the definitions of the diagonal cohomology for the Lie algebra $\on{Vect}M$ with values in $\R$ and $\Om^*(M)$ and the relationship between the cohomology of the double complex $C^*_\De(\on{Vect}M;\Om^*(M)$ and the diagonal cohomology  $H^*_\De(\on{Vect}M)$.

Section \ref{main bicomplex} contains the definition of the diagonal complex 
$\Om^*_\De(D\x M)$, the filtration of $\Om^*_\De(D\x M)$ and the corresponding spectral sequence, and the proof of the cohomology isomorphism 
$H^*(\Om^*_\De(D\x M)^D)=H^*_\De(\on{Vect}M;\Om^*(M)$ (Lemma \ref{iso}). Finally, one obtains the main result on the cohomology $H^*(D)$, when $H^i(M)=0$ for $i>0$ (Corollary \ref{main1}), and on the cohomology $H^*(D_+)$, when $M$ is a closed oriented manifold (Corollary \ref{main2}). 

Throughout the paper $M$ is a connected $n$-dimensional oriented manifold with countable base of $C^\infty$-class, smooth map means $C^\infty$-map and, for a finite or infinite dimensional manifold $X$,  $\Om^*(X)=(\Om^p(X))$ means the de Rham complex of $X$.    
 
\section{Preliminaries}\label{prel}
\subsection{A double complex} 
Let $G$ be a topological group and let $K$ be a left topological $G$-module. Recall (see,for example, \cite{Guic}) that the standard complex $C^*_{\on{cont}}(G,K)=\{C^p_{\on{cont}}(G,K),d\}$ of continuous nonhomogeneous cochains with the differential $d=\{d^p\}$ is defined as follows: for $p>0$, 
$C^p_{\on{cont}}(G,K)$ is the space of continuous maps from $G^p$ to $K$, $C^0(G,K)=K$, and, for $c\in C^p_{\on{cont}}(G,K)$,  we have
\begin{multline*}
(d^pc)(g_1,\dots,g_{p+1})=g_1c(g_2,\dots,g_{p+1})\\+
\sum_{i=1}^p(-1)^{i-1}c(g_1,\dots,g_ig_{i+1},\dots,g_{p+1})+(-1)^{p+1}c(g_1,\dots,g_p),
\end{multline*}
where $g_1,\dots,g_{p+1}\in G$. The $p$\,th cohomology group of this complex is denoted by $H^p_{\on{cont}}(G,K)$ and is called the $p$\,th continuous cohomology group of $G$ with values in $K$.  

Let $G$ be a topological group, $K=\{K^q,\,d^q\}$ a cochain complex  such that each $K^q$ is a left topological $G$-module, and $d^q:K^q\to K^{q+1}$  a $G$-equivariant continuous homomorphism of modules. Later we consider the following construction of a double complex $C^*_{\on{cont}}(G;K)$ and some of its applications (\cite{L2},\cite{L4}). 

Let $d^{p,q}:C^p_{\on{cont}}(G,K^q)\to C^{p+1}_{\on{cont}}(G,K^q)$ be the differential of the standard complex $C^p_{\on{cont}}(G,K^q)$. By the standard way we will consider $C^*_{\on{cont}}(G,K)$ as a  double complex putting $\de_1=\{\de_1^{p,q}\}$, where $\de_1^{p,q}=d^{p,q}$, and defining the second differential $\de_2=\{\de_2^{p,q}\}$ in the following way: for $c\in C^p_{\on{cont}}(G,K^q)$ we 
put 
$$
\de_2^{p,q}c(\cdot)=(-1)^pd^qc(\cdot).
$$
Then $C^*_{\on{cont}}(G,K)=\{C^p_{\on{cont}}(G,K^q)\}$ is a cochain complex with respect to the total differential $\de=\de_1+\de_2$ and the total grading $C^m(G,K)=\otimes_{p+q=m}C^p_{\on{cont}}(G,K^q)\}$. We denote this complex by $C^*_{\on{cont}}(G;K)$ and denote by $H^p_{\on{cont}}(G;K)$ the $p$th cohomology group of this complex. 
 
Let $K^G$ be the subcomplex of $G$-invariant cochains of $K$. Evidently  
$K^G\subset C^0(G,K)$ is a subcomplex of the complex $C^*_{\on{cont}}(G;K)$ and we have the corresponding cohomology homomorphism $H^*(K^G)\to H^*(G;K)$. This cohomology homomorphism  will play an important role in the following constructions. 

Next we mainly consider the case when $K=\{K^q\}$ is a differential graded algebra 
(briefly $DG$-algebra) and the differential $d=\{d^q\}$ is an antiderivation of this algebra of degree $1$. Then $C^*_{\on{cont}}(G;K)$ is a $DG$-algebra also and the total differential is an antiderivation of this algebra of degree $1$. 

Let $f:K'\to K''$ be a $G$-equivariant homomorphism of topological $G$-complexes. It is easy to check that $f$ induces a homomorphism of the corresponding double complexes and, therefore, a homomorphism of complexes $C^*(G;K')\to C^*(G;K'')$.
 
\subsection{The cohomology $H^*(W_n)$}
Let $\mathfrak g$ be a topological Lie algebra and let $E$ be a left topological 
$\mathfrak g$-module. Recall (see, for example, \cite{Guic}) that  the complex 
$C^*_{\on{cont}}(\mathfrak g,E)=\{C^q_{\on{cont}}(\mathfrak g,E),d^q\}$ of standard continuous cochains of $\mathfrak g$ with values in $E$ is defined as follows: $C^q_{\on{cont}}(\mathfrak g,E)$ is the space of continuous skew-symmetric $q$-forms on $\mathfrak g$ with values in $E$ and the differential 
$d^q:C^q_{\on{cont}}(\mathfrak g,E)\to C^{q+1}_{\on{cont}}(\mathfrak g,E)$ is defined by the following formula:
\begin{multline*}
(d^q c)(\xi_1,\dots,\xi_{q+1})=
\sum_{i=1}^{q+1}(-1)^{i-1}\xi_ic(\xi_1,\dots,\widehat{\xi_i},\dots,\xi_{p+1})\\+
\sum_{i,j}(-1)^{i+j}c([\xi_i,\xi_j],\xi_1,\dots,
\widehat{\xi_i},\dots,\widehat{\xi_j},\dots,\xi_{q+1}),
\end{multline*}
where $c\in C^q_{\on{cont}}(\mathfrak g,E)$, $\xi_1,\dots,\xi_{q+1}\in\mathfrak g$, and, as usual, $\hat\xi$ means that the term $\xi$ is omitted.  We denote the cohomology of this complex by $H^*_{\on{cont}}(\mathfrak g,\R)=\{H^p_{\on{cont}}(\mathfrak g,\R)\}$. 

If $E$ is an algebra and the action of $\mathfrak g$ on $E$ is compatible with this algebra structure, the complex 
$C^*_{\on{cont}}(\mathfrak g,E)$ is a graded algebra with respect to the exterior product of forms induced by the product in $E$ and the differential $d=\{d^q\}$ is an antiderivation of the graded algebra $C^*(\mathfrak g,E)$ of degree 1.  

Let $W_n$ be the algebra of formal vector fields in $n$ variables, 
i.e. the topological vector space of $\infty$-jets at $0$ of smooth vector fields on $\R^n$ with the bracket induced by the Lie bracket of vector fields on 
$\R^n$. Consider $\R$ as a trivial $W_n$-module. For brevity, put $C^*(W_n)=C^*_{\on{cont}}(W_n,\R)$ and $H^*(W_n)=H^*_{\on{cont}}(W_n,\R)$.  

Recall some facts about the cohomologies $H^*(W_n)$ and $H^*(W_n,\on{GL}_n(\R)$ 
(\cite{BR}, \cite{F}, \cite{God}). 
Consider the complex $C^*(W_n)=\{C^q(W_n),d^q\}$ of standard continuous cochains of $W_n$ with values in the trivial $W_n$-module $\R$. By definition, $C^*(W_n)$ is a $DG$- algebra and the differential $d$ is an antiderivation of degree 1. For $\xi^i\in\R[[\R^n]]$ and 
$\xi=\sum_{i=1}^n\xi^i\frac{\p}{\p x^i}\in W_n$, put 
$$
c^i_{j_1\dots j_r}(\xi)=\frac{\p^r\xi^i}{\p x^{j_1}\dots\p x^{j_r}}(0),
$$
where $x^i$ $(i=1,\dots,n)$ are the standard coordinates in $\R^n$. By definition, we have  $c^i_{j_1\dots j_r}\in C^1(W_n)$. Moreover, $c^i_{j_1\dots j_r}$  for $r=0,1,\dots$ and $i,j_1\dots j_r=1,\dots,n$ are generators of the $DG$-algebra $C^*(W_n)$.   
Since $d=\{d^q\}$ is an antiderivation of degree 1 of $C^*(W_n)$, it is uniquely determined by the following conditions: 
\begin{equation}\label{dc}
dc^i_{j_1\dots j_r}=
\sum_{o\le k\le r}\sum_{s_1<\dots s_k}\sum_{l=1}^nc^i_{lj_1\dots\widehat{j_{s_1}}\dots
\widehat{j_{s_k}}\dots j_r}\wedge c^l_{j_{s_1}\dots j_{s_k}}.
\end{equation}

The group $\on{GL}_n(\R)$ acts naturally on $C^*(W_n)$ and its Lie algebra 
$\mathfrak{gl}_n(\R)$ as the Lie algebra of vector fields with linear components is a subalgebra of $W_n$. Then we have the complex $C^*(W_n,\on{GL}_n(\R))$ of relative cochains of the Lie algebra $W_n$ with respect to $\on{GL}_n(\R)$ and the cohomology $H^*(W_n,\on{GL}_n(\R))$ of this complex. Similarly, we have the complex $C^*(W_n,\on{O}(n))$ of relative cochains of the Lie algebra $W_n$ with respect to the orthogonal group $\on{O}(n)\subset\on{GL}_n(\R)$ and the cohomology $H^*(W_n,\on{O}(n))$. 

Put 
$$
\ga=(c^i_j),\quad \Psi^i_j=\sum_{k=1}^nc^i_{jk}\wedge c^k,\quad\on{and}\quad  \Psi=(\Psi^i_j).
$$  
It is known that 
$$
\Psi_p=\on{tr}(\underbrace{\Psi\wedge\dots\wedge\Psi}_{\text {p times}})\quad (p=1,\dots,n).
$$
are cocycles of $C^*(W_n,\on{GL}_n(\R))$ and the cohomology classes of these cocycles  generate $H^*(W_n,\on{GL}_n(\R))$. The cohomology class of $\Psi_k$ is called $k$th formal  Pontrjagin class. 

Put
$$
\ga_p=\on{tr}(\underbrace{\ga\wedge\dots\wedge\ga}_{\text {2p-1 times}})\quad (p=1,\dots,n).
$$
By definition, $\ga_p\in C^{2p-1}(W_n)$. Consider the inclusion $\mathfrak{gl}_n(\R)\subset W_n$. It is known that $\ga_p$ as a cochain of the complex $C^*(\mathfrak{gl}_n(\R),\R)$ is a cocycle and   
the ring $H^*(\mathfrak{gl}_n(\R),\R)$ is the exterior algebra of its subspace spanned by the cohomology classes $\ga_p$ for $p=1,\dots,n$. Moreover, there is a $(2p-1)$-cochain $\Ga_p$ of the complex $C^*(W_n)$ such that the restriction of $\Ga_p$ to $\mathfrak{gl}_n(\R)$ equals $\ga_p$ and $d\Ga_p=\Psi_p$. Consider the $DG$-subalgebra of $DG$-algebra  $C^*(W_n)$ generated by $\Ga_p$ and $\Psi_p$ for $p=1,\dots,n$. Then the inclusion of this subalgebra into $C^*(W_n)$ induces an isomorphism of the cohomologies. Moreover, the cohomology classes of the cocycles  
\begin{equation}\label{basis}
\Ga_{p_1}\wedge\dots\wedge\Ga_{p_l}\wedge\Psi_{r_1}\wedge\dots\wedge\Psi_{r_m}
\end{equation} 
for 
$1\le p_1<\dots<p_l\le n$, $1\le r_1\le\dots\le r_m\le n$, $p_1\le r_1$, $r_1+\dots+r_m\le n$, and $p_1+r_1+\dots+r_m>n$ give a basis of $H^*(W_n)$ (the so-called Vey basis) as a vector space. This implies that $H^m(W_n)=0$ whenever $m<2n+1$ or $m>n(n+2)$.

\subsection{The space of frames of infinite order and the Gelfand-Kazhdan 
form }\label{S(M)} 
Let $M$ be a connected orientable $n$-dimensional smooth manifold. Denote by $S(M)$ the space of  frames of infinite order of $M$, i.e. $\infty$-jets at $0$ of germs at $0$ of diffeomorphisms from $\R^n$ into $M$. It is known that $S(M)$ is a manifold with model space $\R^\infty$ (\cite{BR}). Recall that we denote by $D$ the group of diffeomorphisms $\on{Diff}M$ of $M$. We put, for $g_1,g_2\in D$, $g_1g_2=g_2\o g_1$. Then the standard action of $D$ on $M$ is a right action. Evidently this action of $D$ is naturally extended to an action of $D$ on $S(M)$. 

Define the canonical Gelfand-Kazhdan 1-form $\om$ with values in $W_n$ on $S(M)$ 
(\cite{G-K} and \cite{BR}). Let $\tau$ be a tangent vector at $s\in S(M)$ and  let $s(t)$ 
be a path on $S(M)$ such that $\tau=\frac{ds}{dt}(0)$. One can represent $s(t)$ by a smooth family $k_t$ of germs at $0$ of diffeomorphisms $\R^n\to M$, i.e. $s(t)=j^\infty_0k_t$. Then put
$$
\om(\tau)=-j_0^\infty\frac{d}{dt}(k_0^{-1}\o k_t)(0).
$$
  
\begin{theorem}\label{G-K}(\cite{BR}, \cite{G-K}, \cite{GS}) The form $\om$ satisfies the following conditions: 
\begin{enumerate}
\item $\om$ induces a topological isomorphism between the tangent space $T_s$ of $S(M)$ at $s\in S(M)$ and $W_n$;  
\item $d\om=-\frac12[\om,\om]$\quad (the Maurer-Cartan condition);
\item The form $\om$ is $D$-invariant. 
\end{enumerate}
\end{theorem}
This theorem and the theorem on the covering isotopy of imbeddings of disks (\cite{Cerf},  \cite{Pale}) implies the following 
\begin{corollary}\label{disk}
The group $D$ acts transitively on $S(M)$.
\end{corollary}
  
Let $c\in C^q(W_n)$. 
For each $s\in S(M)$ and $X_1,\dots,X_q\in T_s$, put 
$$
\om_c(X_1,\dots,X_q)=c(\om(X_1),\dots,\om(X_q)).
$$ 
By theorem \ref{G-K}, $c\mapsto\om_c$ is an injective homomorphism of the complexes $\al:C(W_n)\to\Om^*(S(M))$. 
By theorem \ref{G-K} and corollary \ref{disk}, there is a one-to-one correspondence between the space of $D$-invariant forms on $S(M)$ and the space of continuous skew-symmetric forms on the tangent space $T_s$ at $s\in S(M)$. Then $\al(C(W_n))=\Om^*(S(M))^D$, where 
$\Om^*(S(M))^D$ is the subcomplex  of $D$-invariant forms from $\Om^*(S(M))$. 
Moreover, we have $\al(C(W_n,\on{O}(n)))=\Om^*(S(M)/\on{O}(n))^D$. 
It is easy to check that $(\be^i_j)=-(\al(c^i_j))$ is a connection form on a principal $\on{GL}_n(\R)$-bundle $S(M)\to S(M)/\on{GL}_n(\R)$. 

Consider a Riemannian metric on $M$ and the corresponding Levi-Civita connection. 
Denote by $O(M)$ the principal $\on{O}(n)$-bundle  of orthogonal tangent frames on $M$. For each frame $r\in O(M)$ at $x\in M$, denote by $\si(r)$ the $\infty$-jet at $x$ of the inverse of the geodesic chart with center at $x$ uniquely determined by $r$. Then we have a smooth map $\si:O(M)\to S(M)$. 

Denote by $\theta=(\theta^i)$ the canonical 1-form $\theta=(\theta^i)$ with values in $\R^n$ on $O(M)$ and by $\theta=(\theta^i_j)$ the form of the Levi-Civita connection.
It is easy to check that, for $\om^i_j=\al(c^i_j)$ and $\om^i_j=\al(c^i_j)$, we have  \begin{equation}\label{con}
\si^*\om^i=-\theta^i\quad \text{and}\quad \si^*\om^i_j=-\theta^i_j.
\end{equation}

Let $c^i_{j_1\dots j_r}\in C^1(W_n)$. For a vector field $X$ on $M$, denote by $\tilde X$
the corresponding extension of $X$ to $S(X)$. Put 
$\om^i_{j_1\dots j_r}=\al(c^i_{j_1\dots j_r})$. 
By definition, we have
\begin{equation}\label{om(X)} 
\om^i_{j_1\dots j_r}(\tilde X(s))=-\frac{\p^r X^ i}{\p x^{j_1}\dots\p x^{j_r}}(0),
\end{equation}
where $x\in M$ and the right hand side is calculated in the coordinates $x^i$ presenting $s\in S(M)$. 

\subsection{The cohomology $H^*(W_n,\on{O}(n))$}

Consider the basic cocycles \eqref{basis}. For odd $p$ it is possible to choose $\ga_p$ and $\Ga_p$ so that the restriction of $\ga$ is a cocycle of $C^*(W_n,\on{O}(n))$ and 
$\Ga_p\in C^*(W_n,\on{O}(n))$ (\cite{F}, \cite{God}). In this case we replace such $\ga_p$ and $\Ga_p$ by $\la_p$ and $\La_p$ respectively. 

Consider the $DG$-subalgebra of $DG$-algebra  $C^*(W_n,\on{O}(n))$ generated by $\La_p$ for odd $p\le n$ and $\Psi_p$ for $p=1,\dots,n$. Then the inclusion of this subalgebra into $C^*(W_n,\on{O}(n))$ induces an isomorphism of the cohomologies. Besides, the cohomology classes of cocycles 
\begin{equation}\label{basis1}
c_{p_1\dots p_l,r_1\dots r_k}
=\La_{p_1}\wedge\dots\wedge\La_{p_l}\wedge\Psi_{r_1}\wedge\dots\wedge\Psi_{r_k}
\end{equation} 
for odd $p_i$ such that 
$1\le p_1<\dots<p_l\le n$, $1\le r_1\le\dots\le r_k\le n$, $p_1\ge r_1$, $r_1+\dots+r_k\le n$, and $p_1+r_1+\dots+r_k>n$ give a basis of the cohomology $H^*(W_n,\on{O}(n))$ as a vector space in dimensions $>2n$. In particular, for each dimension $m>2n$ we have 
$H^m(W_n,\on{O}(n))\subset H^m(W_n)$. Moreover, for $m>n$ the cohomology group $H^m(W_n,\on{O}(n))$ may be nontrivial only if $2n+1\le m\le\frac{n(n+3)}{2}$, when $n$ is even, and $2n+1\le m\le\frac{n(n+5)}{2}$, when $n$ is odd.       

Next we find the explicit expressions for $\La_p$ for odd $p$.

Put $\la^i_j=\frac12(c^i_j+c^j_i)$, $\la=(\la^i_j)$, and 
$$
\la_p=\on{tr}(\underbrace{\la\wedge\dots\wedge\la}_{\text {2p-1 times}})\quad (p=1,\dots,n).
$$   
It is easy to check that the restriction of $\la_p$ to $\mathfrak{gl}_n(\R)$ belongs to  $C^*(\mathfrak{gl}_n(\R),\on{O}(n))$ and $\la_p=0$ for even $p$. 

Let $P$ be a homogeneous $\on{GL}_n(\R)$-invariant polynomial of degree $p$ on the Lie algebra 
$\mathfrak{gl}_n(\R)$. We denote by $\bar P$ the symmetric $p$-linear form on 
$\mathfrak{gl}_n(\R)$ such that $\bar P(X,\dots,X)=P(X)$. Then for $X,X_1,\dots,X_p\in\mathfrak{gl}_n(\R)$ we have 
\begin{equation}\label{Ginv}
\sum_{i=1}^p\bar P(X_1,\dots,[X,X_i],\dots,X_p)=0. 
\end{equation}   
If $\om$ is a linear form with values in $\mathfrak{gl}_n(\R)$ and, for $i=1,\dots,p$, $\omega_i$ is a $k_i$-linear form with values in $\mathfrak{gl}_n(\R)$, \eqref{Ginv} implies the equality
\begin{equation}\label{Pinv}
\sum_{i=1}^p(-1)^{k_1+\dots+k_{i-1}+1}
P(\om_1,\dots,[\om,\om_i],\dots,
\om_p)=0. 
\end{equation}

Put $\al^i_j=\frac12(c^i_j-c^j_i)$ and $\al=(\al^i_j)$.
By definition, $\la^i_j,\al^i_j\in C^1(W_n)$. By \eqref{dc}, we have
$$
d\la^i_j=\sum_{k=1}^n\left(\la^i_k\wedge\psi^k_j+\psi^i_k\wedge\la^k_j\right)+
\frac12(\Psi^i_j+\Psi^j_i).
$$
This formula could be rewritten as follows: 
\begin{equation}\label{dla}
d\la=[\al,\la]+\frac12(\Psi+\Psi^t).
\end{equation}
Similarly, we have 
\begin{gather}
d\Psi=[\la,\Om]+[\al,\Psi],\label{dOm}\\
d\Psi^t=-[\la,\Psi^t]+[\al,\Psi^t].\label{dOmt}
\end{gather}

Let $P$ be a $\mathfrak{gl}_n(\R)$-invariant $(a+b+c)$-linear form  on $\mathfrak{gl}_n(\R)$. Then 
$$
P(\underbrace{\la,\dots,\la}_{\text{$a$ times}},\underbrace{\Psi,\dots,
\Psi}_{\text{$b$ times}},\underbrace{\Psi^t,\dots,\Psi^t}_
{\text{$c$\ times}}) 
$$
is a cochain of $C^{a+2b+2c}(W_n,\on{O}(n))$. Since $[\la[\la,\la]]=0$, by \eqref{Pinv}, \eqref{dla}, \eqref{dOm}, and \eqref{dOmt}, we have
\begin{equation}\label{dP}
\begin{aligned}
&dP(\la,\dots,\la,\Psi,\dots,\Psi,\Psi^t,\dots,\Psi^t) \\
&=\sum_{i=1}^a(-1)^{i-1}P(\la,\dots,\Psi+\Psi^t,\dots,\la,\Psi,\dots,
\Psi,\Psi^t,\dots,\Psi^t) \\
&+(-1)^a\sum_{j=1}^bP(\la,\dots,\la,\Psi,\dots,[\la,\Psi],\dots,
\Psi,\Psi^t,\dots,\Psi^t) \\
&-(-1)^a\sum_{k=1}^cP(\la,\dots,\la,\Psi,
\dots,\Psi,\Psi^t,\dots,[\la,\Psi^t],\dots,\Psi^t).
\end{aligned} 
\end{equation}

Let $Q$ be a homogeneous invariant polynomial of degree $p$ on $\mathfrak{gl}_n(\R)$. For  brevity, put $Q(X,X')=\bar Q(X,X',\dots,X')$ and 
$Q(X,X',X'')=\bar Q(X,X',X'',\dots,X'')$.  
Let $Q_p(X)=\on{tr}X^p$. Then we have 
$$
\bar Q_p(X_1,\dots,X_p)=\frac1{p\,!}\sum_{\si\in S_p}\on{tr}\left(X_{\si(1)}\dots X_{\si(p)}\right),
$$
where $S_p$ is the symmetric group. 

For odd $p$, put $\Psi(t)=\frac12t\Psi+\frac12(t-1)\Psi^t+(t-t^2)[\la,\la]$. It is evident that 
$$
\int_0^1\frac{d}{dt}Q_p(\Psi(t)\,dt=Q_p\Bigl(\frac12\Psi\Bigr)+Q_p\Bigl(\frac12\Psi^t\Bigr)=\frac1{2^{p-1}}\Psi_p.
$$
On the other hand, by \eqref{dla}, \eqref{dOm}, \eqref{dOmt}, \eqref{Pinv}, and  \eqref{dP} we have
\begin{multline*} 
\int_0^1\frac{d}{dt}Q_p(\Psi(t))\,dt =
p\int_0^1\left(Q_p(2^{-1}(\Psi+\Psi^t),\Psi(t))+(1-2t)Q_p([\la,\la],\Psi(t))\right)\,dt\\=
p\int_0^1\Bigl(Q_p(2^{-1}(\Psi+\Psi^t),\Psi(t))-2^{-1}(p-1)t(2t-1)Q_p(\la,[\la,\Psi],\Psi(t))\\-(p-1)(t-1)2^{-1}(2t-1)Q_p(\la,[\la,\Psi^t],\Psi(t))\Bigr)\,dt=
p\,d\int_0^1Q_p(\la,\Psi(t))\,dt.
\end{multline*}

It is easy to see that the restriction of $\int_0^1Q_p(\la,\Psi(t))\,dt$ to $C^*(\mathfrak{gl}_n(\R),\R))$ equals $\frac{p\,!}{(p+1)\dots(2p+1)}\la_p$. Therefore, one can put 
\begin{equation}\label{Lap}
\La_p=\frac{p\,p\,!\,2^{p-1}}{(p+1)\dots(2p+1)}\int_0^1Q_p(\la,\Psi(t))\,dt.
\end{equation}

\section{Cocycles on diffeomorphism groups}\label{Dcocycles}

\subsection{The definition of cocycles}\label{cocycles}

Later we need the following theorem.    
\begin{theorem}\label{split}
Let $M$ be a connected orientable manifold. Then, in the category of topological vector spaces, 
for each $p>0$ we have the following decomposition 
$$
\Om^p(M)=d\Om^{p-1}(M)\oplus H^p(M)\oplus\Om^p(M)/Z^p(M),
$$
where $Z^p(M)$ is the space of closed $p$-forms. If $H^p(M)=0$,
$d\Om^{p-1}(M)=Z^p(M)$ and $\Om^p(M)/Z^p(M)$ are Frech\'et spaces.
\end{theorem}

\begin{proof} For a compact $M$ the statement follows from the Hodge decomposition
for the identity operator 1 on $\Om^p(M)$: 
$1=d\o\de\o G\oplus H^p(M)\oplus\de\o d\o G$ (see, for example, \cite{Rh}).

For a noncompact $M$, the statement follows from Palamodov's theorem (\cite{Pal},
Proposition 5.4).
\end{proof}

The action of $D$ on $M$ induces the natural right actions of $D$ on $S(M)$ and   $S(M)/\on{O}(n)$ and the induced left actions on the de Rham complexes $\Om^*(S(M))$ and  $\Om^*(S(M)/\on{O}(n))$. Consider the $D$-module $\Om^*(S(M)/\on{O}(n))$ and the corresponding complex $C^*_{\on{cont}}(D;\Om^*(S(M)/\on{O}(n)))$ . 
If the manifold $M$ is oriented, denote by $D_+$ the subgroup of $D$ of diffeomorphisms preserving the orientation of $M$. 

Further we consider the complex $\Om^*(S(M)/\on{O}(n))^D$ as a subcomplex of the complex $C^*(D;\Om^*(S(M)/\on{O}(n))$. 

The proof of the next theorem shows how to construct cocycles on the group $D$ from the cocycles of the complex $\Om^*(S(M)/\on{O}(n))^D$. 

Consider a structure of $D$-module on a cohomology group $H^i(M)$ induced by the action of the group $D$ on $\Om^*(M)$. 
\begin{theorem}\label{group cocycle}
Let $\nu$ be an $m$-cocycle of the complex $\Om^*(S(M)/\on{O}(n))^D$ and $0\le p\le m-1$. Assume that $H^{m-p}(M)=\dots =H^m(M)=0$  and $H^{m-p-1}(M)\ne 0$. Then the following statements are true:
\begin{enumerate}
\item There are cochains 
$\nu_{i,m-i-1}\in C^i_{\on{cont}}(D;\Om^{m-i-1}(S(M)/\on{O}(n))$ $(i=0,\dots,p)$ 
such that $\nu+\de(\nu_{0,m-1}+\dots+\nu_{p,m-p-1})=\de_1\nu_{p,m-p-1}$;
\item $c(\nu)=\de_1\nu_{p,m-p-1}\in C^{p+1}_{\on{cont}}(D;\Om^{m-p-1}(S(M)/\on{O}(n))$ defines a $(p+1)$-cocycle of the complex $C^{p+1}_{\on{cont}}(D,H^{m-p-1}(M)$.  
The cohomology class $h(\nu)$ of the cocycle $c(\nu)$ depends only on the cohomology class of the cocycle $\nu$ in the complex $\Om^*(S(M)/\on{O}(n))^D$;
\item The map $\nu\mapsto h(\nu)$ induces a linear map 
$$
H^m(\Om^*(S(M)/\on{O}(n))^D)\to H^{p+1}_{\on{cont}}(D,H^{m-p-1}(M)).
$$ 
\end{enumerate}
\end{theorem} 
\begin{proof}
Denote by $R(M)$ the space of $\infty$-jets of of germs of Riemannian metrics at points of $M$. Let $g_x$ be the $\infty$-jet at $x\in M$ of a germ at $x$ of Riemannian metric $g$. For an orthogonal frame $r$ at $x$, denote by $s(r)$ the $\infty$-jet at the center $x$ of the inverse of the geodesic chart defined by $r$ and the Riemannian metric $g$. It is evident that $s(r)$ is uniquely determined by the $\infty$-jet $g_x$. Then the map $r\mapsto s(r)$ induces a diffeomorphism $R(M)\to S(M)/\on{O}(n)$. Next we identify $R(M)$ and $S(M)/\on{O}(n)$ by this diffeomorphism. 

Consider some Riemannian metric $g_0$ on $M$ and the smooth map $\si:O(M)\to S(M)$ defined by $g_0$ in \ref{S(M)}. Let $\tilde\si:M\to S(M)/\on{O}(n)$ be the smooth map induced by $\si$.  Consider a smooth map $F:[0,1]\x R(M)\to R(M)$ defined by 
$$
(t,g_x)\mapsto (tg_0+(1-t)g)_x=tg_{0,x}+(1-t)g_x,
$$ 
where $g_{0,x}$ and $g_x$ are the $\infty$-jets at $x$ of germs at $x$ of Riemannian metrics $g_0$ and $g$. It is clear that $F$ is a smooth homotopy between the identity map 
of $R(M)$ and the composition of the projection $p:R(M)\to M$ and the section $\tilde\si$. 
Then we have the standard homotopy operator 
$\Ph^p:\Om^p(S(M)/\on{O}(n))\to \Om^{p-1}(S(M)/\on{O}(n))$ given by 
$\om\mapsto\int_0^1\on{i}(\frac{\p}{\p t})F^*\om ~dt$, where $\on{i}(X)$ denote the inner product by a vector field $X$. In particular, we have 
$H^*(S(M)/\on{O}(n))=H^*(M)$. 

First we indicate some canonical construction of the sequence of cochains 
$$
\nu_{i,m-i-1}\in C^i_{\on{cont}}(D;\Om^{m-i-1}(S(M)/\on{O}(n))\quad (i=0,\dots,p)
$$ 
satisfying the conditions of the theorem.
 
By assumption, we have $\Ph^m(M)=0$ and $p\ge 0$.  
Suppose that $m>n$. Then, for the $(m-1)$-form $\nu_{0,m-1}=-H^m(\nu)$, we have $\nu=-\de_2\nu_{0,m-1}$. If $m\le n$, we use theorem \ref{split} to get the $(m-1)$-form $\nu_{0,m-1}$ satisfying the same equality $\nu=-\de_2\nu_{0,m-1}$. In the both cases one can assume that 
$\nu_{0,m-1}\in C^0_{\on{cont}}(D;\Om^{m-1}(M))$. Then we have
$\nu+\de\nu_{0,m-1}=\de_1\nu_{0,m-1}$ and $\de_2\de_1\nu_{0,m-1}=-\de_1\de_2\nu_{0,m-1}=\de_1\nu=0$. 

Let $p\ge 1$ and, therefore, $H^{m-1}(M)=0$.  
Suppose that $m-1>n$. Then, for the cochain $\nu_{1,m-2}=\Ph^{m-1}(\nu_{1,m-1})$, we have 
$\de_1\nu_{0,m-1}=-\de_2\nu_{1,m-2}$.   
If $m-1\le n$, we use theorem \ref{split} to get a cochain $\nu_{1,m-2}\in C^1_{\on{cont}}(D;\Om^{m-2}(M))$ satisfying the same equality $\de_1\nu_{0,m-1}=-\de_2\nu_{1,m-2}$. 
Thus, we have
$\nu+\de(\nu_{0,m-1}+\nu_{1,m-2})=\de_1\nu_{1,m-2}$ and $$
\de_2\de_1\nu_{1,m-2}=-\de_1\de_2\nu_{1,m-2}=\de_1^2\nu_{0,m-1}=0.
$$ 

Using the conditions
$$
H^{m-2}(S(M)/\on{O}(n))=\dots =H^{m-p}(S(M)/\on{O}(n))=0
$$
and proceeding in the same way we get
for $i=1,\dots,p$  the cochains $\nu_{i,m-i-1}\in C^i_{\on{cont}}(D;\Om^{m-i-1}(S(M)/\on{O}(n))$
such that
\begin{equation}\label{sequence}
\de_1\nu_{i-1,m-i}+\de_2\nu_{i,m-i-1}=0
\end{equation}
and so
$$
\nu+\de(\nu_{0,m-1}+\dots+\nu_{p,m-p-1})=\de_1\nu_{p,m-p-1}\in C^{p+1}(D,\Om^{m-p-1}(S(M)/\on{O}(n)).
$$
Moreover, we have
$$
\de_2\de_1\nu_{p,m-p-1}=-\de_1\de_2\nu_{p,m-p-1}=\de_1^2\nu_{p-1,m-p}=0.
$$
Consider $H^{m-p-1}(S(M)/\on{O}(n))=H^{m-p-1}(M)$ as a $D$-module with respect to the natural action of $D$.
Then the cochain $\de_1\nu_{p,m-p-1}$ defines a $(p+1)$-cocycle $c(\nu)$ on $D$ with values in $H^{m-p-1}(M)$. The cohomology class of $c(\nu)$ is denoted by $h(\nu)$.
We claim that the cohomology class $h(\nu)$ depends only on the
cohomology class of $\nu$ in the complex $\Om^*(S(M)/\on{O}(n))^D$.

If we replace the form $\nu$ by a form $\nu+d\nu_1$, where
$\nu_1\in\Om^{m-1}(S(M)/\on{O}(n))\cap\Om^*(S(M)/\on{O}(n))^D$, one can replace the sequence
$\nu_{0,m-1},\dots,\nu_{p,m-p-1}$ by the sequence
$\nu_{0,m-1}-\nu_1,\nu_{1,m-2},\dots,\nu_{p,m-p-1}$ and obtain the same cochain
$\nu_{p,m-p-1}$ at the end.

Consider another sequence $\bar\nu_{0,m-1},\dots,\bar\nu_{p,m-p-1}$
$(i=0,\dots,p)$ such that 
$$
\nu=-\de_2\bar\nu_{0,m-1}\quad {\text and}\quad 
\de_1\bar\nu_{i-1,m-i}+\de_2\bar\nu_{i,m-i-1}=0
$$ 
for $i=1,\dots,p$.
The same arguments as above show that  
$$
\bar\nu_{0,m-1}=\nu_{0,m-1}+\de_2\si_{0,m-2},
$$
where $\si_{0,m-2}\in C^0(D;\Om^{m-2}(M))$. If $p=1$, we have
$\de_1\bar \nu_{0,m-1}=\de_1\nu_{0,m-1}-\de_2\de_1\si_{0,m-2}$ and we are done.
If $p>1$ we have
$$
\de_1\bar\nu_{0,m-1}=\de_1\nu_{0,m-1}-\de_2\de_1\si_{0,m-2}=-\de_2(\nu_{1,m-2}
+\de_1\si_{0,m-2})=-\de_2\bar\nu_{1,m-2}.
$$

For $i=1,\dots,p-1$ proceeding in the same way we get the cochains
$\si_{i,m-i-2}\in C^i_{\on{cont}}(D;\Om^{m-i-2}(S(M)/\on{O}(n))$ such that
$$
\bar\nu_{i,m-i-1}=\nu_{i,m-i-1}+\de_1\si_{i-1,m-i-1}+\de_2\si_{i,m-i-2}.
$$
In particular, we have
$$
\bar\nu_{p,m-p-1}=\nu_{p,m-p-1}+\de_1\si_{p-1,m-p-1}+\de_2\si_{p,m-p-2}
$$
and 
$$\de_1\bar\nu_{p,m-p-1}=\de_1\nu_{p,m-p-1}-\de_2\si_{p,m-p-2}.
$$ 
Thus, the cocycles $\de_1\bar\nu_{p,m-p-1}$ and $\de_1\nu_{p,m-p-1}$ define the same
cohomology class of $H^{p+1}_{\on{cont}}(D,H^{m-p-1}(M))$. 

The last statement of the theorem follows from the construction. 
\end{proof}
Theorem \ref{group cocycle} implies the following corollaries.
\begin{corollary}\label{cocycle1}
Let $H^i(M)=0$ for $i>0$. Then, for each nontrivial $m$-cocycle $\nu$ of the complex $C^*(W_n,\on{O}(n))=\Om^*(S(M)/\on{O}(n))^D$, the cocycle $c(\nu)$ is an $m$-cocycle of the complex $C^*(D)$. The map $\nu\mapsto h(\nu)$ induces a linear map 
$$
H^m(\Om^*(S(M)/\on{O}(n))^D)\to H^m(D).
$$ 
\end{corollary}

\begin{corollary}\label{cocycle2}
Let $M$ be a closed oriented $n$-dimensional manifold and let $D_+$ be the group of diffeomorphisms of $M$ preserving the orientation of $M$. For $m>n$ and each  $m$-cocycle $\nu$ of the complex $\Om^*(S(M)/\on{O}(n))^D$, $\int_Mc(\nu)$ is a $(m-n)$-cocycle of the complex $C^*(D_+)$ presenting the cohomology class $h(\nu)$. The map $\nu\mapsto h(\nu)$ induces a linear map $H^m(\Om^*(S(M)/\on{O}(n))^D)\to H^{m-n}(D_+)$.
\end{corollary}
\begin{proof} Since $H^i(M)=0$ for $i>n$ and $H^n(M)=\R$ one could apply theorem \ref{group cocycle}. Evidently, all statements of theorem \ref{group cocycle} are true if we replace the group $D$ by the group $D_+$. But then $\int_Mc(\nu)$ presents the cohomology class $h(\nu)$. 
\end{proof}

The main aim of the paper is to prove that, for $m>2n$, the linear maps 
$$
H^m(W_n,\on{O}(n))=H^m(\Om^*(S(M)/\on{O}(n))^D)\to H^m(D)
$$
and 
$$
H^m(W_n,\on{O}(n))=H^m(\Om^*(S(M)/\on{O}(n))^D)\to H^{m-n}(D_+)
$$
induced by the maps $\nu\to h(\nu)$ and $\nu\to \int_Mc(\nu)$ from corollaries \ref{cocycle1} and \ref{cocycle2} are monomorphisms.

\subsection{The diagonal cohomology of the Lie algebra $\on{Vect}M$}

Denote by $\on{Vect}M$ the Lie algebra of smooth vector fields on $M$. We consider $\on{Vect}M$ as a topological vector space with respect to the $C^\infty$-topology. 
Consider the complex $C^*(\on{Vect}M)$ of standard continuous cochains of $\on{Vect}M$ with values in the trivial $\on{Vect}M$-module $\R$ and its subcomplex formed by chains 
$c\in C^*(\on{Vect}M)$ such that for $X_1,\dots,X_q\in\on{Vect}M$ we have $c(X_1,\dots,X_q)=0$ whenever $\cap_{i=1}^q\on{supp}X_i=\emptyset$. This subcomplex is denoted by $C^*_\De(\on{Vect}M)$ and is called the diagonal subcomplex of $C^*(\on{Vect}M)$ (\cite{F}). The cohomology of $C^*_\De(\on{Vect}M)$ is denoted by $H^*_\De(\on{Vect}M)$ and is called the diagonal cohomology of $\on{Vect}M$ with values in $\R$. 

Consider the left action of the Lie algebra $\on{Vect}M$ on the de Rham complex $\Om^*(M)$ by the Lie derivatives: $\om\mapsto\on{L}_X\om$, where $\om\in\Om^*(M)$, $X\in\on{Vect}M$, and $\on{L}_X$ is the Lie derivative with respect to the vector field $X$. Let $C^*(\on{Vect}M,\Om^*(M))$ be the complex of standard cochains of the topological Lie algebra $\on{Vect}M$ with values in the $\on{Vect}M$-module $\Om^*(M)$. 

A cochain $c\in C^p(\on{Vect}M,\Om^q(M))$ is called diagonal if, for $X_1,\dots,X_p\in\on{Vect}M$, the value $c(X_1,\dots,X_p)$ at $x\in M$ depends only on  the germs of $X_1,\dots,X_p$ at $x$. 
By \cite{P}, this condition is equivalent to the following one: for $X_1,\dots,X_p\in\on{Vect}M$ and $x\in M$, the value $c(X_1,\dots,X_p)$ at $x\in M$ depends only on $\infty$-jets of $X_1,\dots,X_p$ at $x$. 
Denote by $C^*_\De(\on{Vect}M,\Om^*(M))$ the set of diagonal cochains of the complex $C^*(\on{Vect}M,\Om^*(M))$. It is easy to see that $C^*_\De(\on{Vect}M,\Om^*(M))$ is a subcomplex of the complex $C^*(\on{Vect}M,\Om^*(M))$ The subcomplex $C^*_\De(\on{Vect}M,\Om^*(M))$ is called the complex of diagonal cochains of the topological Lie algebra $\on{Vect}M$ with values in the $\on{Vect}M$-module $\Om^*(M)$. 
The cohomology of $C^*_\De(\on{Vect}M,\Om^*(M))$ is denoted by $H^*_\De(\on{Vect}M,\Om^*(M))$ and is called the diagonal cohomology of $\on{Vect}M$ with values in $\Om^*(M)$.  

Denote by $\on{Vect}_cM$ the Lie algebra of vector fields on $M$ with compact supports.  
By definition, each diagonal cochain cochain $c\in C^p(\on{Vect}M,\Om^q(M))$ is uniquely determined by its values on the subalgebra $\on{Vect}_cM$ of the Lie algebra $\on{Vect}M$.  
Define the complex $C^*_\De(\on{Vect}_cM,\Om^*(M))$ of diagonal cochains of the topological Lie algebra $\on{Vect}_cM$ with values in the $\on{Vect}_cM$-module $\Om^*(M)$ similarly.  Remark that the map $C^*_\De(\on{Vect}M,\Om^*(M))\to C^*_\De(\on{Vect}_cM,\Om^*(M))$ induced by the inclusion $\on{Vect}_c(M)\subset\on{Vect}(M)$ is an isomorphism of complexes. 

Since the Lie derivative and the exterior derivative on $\Om^*(M)$ commute, one can endow 
$C_\De^*(\on{Vect}M,\Om^*(M))$ with the second differential induced by the exterior derivative on $\Om^*(M)$. We denote by $C^*_\De(\on{Vect}M;\Om^*(M))$ the corresponding double complex with respect to the total differential and denote by $H^*_\De(\on{Vect}M;\Om^*(M))$ the cohomology of this complex. By the remark above, we see that the inclusion $\on{Vect}_c(M)\subset\on{Vect}(M)$ induces an isomorphism of double complexes 
$$
C^*_\De(\on{Vect}_cM;\Om^*(M))\to C^*_\De(\on{Vect}M;\Om^*(M)).
$$
 It is clear that $C^*_\De(\on{Vect}M;\Om^*(M))$ is a $DG$-algebra and its differential is an antiderivation of degree $1$.  

Let $M$ be a closed oriented $n$-dimensional manifold. Consider a map 
$$
\psi:C^*_\De(\on{Vect}M;\Om^*(M))\to C^*(\on{Vect}M)
$$ 
defined as follows:
for $0\le q <n$ put $\psi(C^*_\De(\on{Vect}M;\Om^q(M)))=0$, for $q=n$ and 
$c\in C^*_\De(\on{Vect}M;\Om^n(M))$, put 
$\psi(c)=\int_Mc$.  
 
\begin{theorem}\label{L-G}(\cite{Guil},\cite{L1})
The map $\psi$ is a homomorphism of complexes which induces an isomorphism $H^{n+p}_\De(\on{Vect}M;\Om^*(M))=H^p_\De(\on{Vect}M)$ for any $p\ge 0$.
\end{theorem}

\section{The main double complex}\label{main bicomplex}

\subsection{The complex $\Om^*(D\x M)$ and its diagonal subcomplex}

By \cite{KM}, $\on{Vect}_cM$ is the 
tangent space $T_e(D)$ of $D$ at the neutral element $e\in D$. Then, for $X\in\on{Vect}_cM$ and $g\in D$, the corresponding left invariant vector field $X^l$ on $D$ is the map $g\mapsto (L_g)_*X$ and the corresponding right invariant vector field $X^r$ on $D$ is the map $g\mapsto (R_g)_*X$, where $L_g$ is a left translation and $R_g$ is a right translation on $D$. Note that, by our assumption, our multiplication in $D$ differs from the multiplication in \cite{KM} by the order of factors.  
By \cite{KM}, for $X,Y\in \on{Vect}_cM$ we have $[X^r,Y^r]=-[X,Y]^r$.

Consider the infinite-dimensional manifold $D\x M$. Then we may consider the Lie algebra 
$V(D\x M)=\on{Vect}_c^rM\oplus\on{Vect}M$  
as a subalgebra of the Lie algebra of vector fields on $D\x M$. 

Put $\widetilde{\on{Vect}}_cM=\{-X^r+X; X\in\on{Vect}_cM\}$ and $\tilde X=-X^r+X$.       
Then $\widetilde{\on{Vect}}_cM$ is a subalgebra of the Lie algebra of vector fields on 
$D\x M$, which is isomorphic to the Lie algebra $\on{Vect}_cM$. Evidently, we have 
$V(D\x M)=\widetilde{\on{Vect}}_cM\oplus\on{Vect}M$.  

Consider the de Rham complex $\Om^*(D\x M)$.  
It is clear that each $m$-form $\om\in\Om^*(D\x M)$ is uniquely determined by the corresponding continuous $m$-linear skew-symmetric form on $V(D\x M)$ with values in the ring of smooth functions on $D\x M$. 

Introduce a bigrading of $\Om^*(D\x M)=\oplus_{p,q}\Om^{pq}(D\x M)$ induced by the decomposition $V(D\x M)=\widetilde{\on{Vect}}_cM\oplus\on{Vect}M$. 
Then $\Om^{pq}(D\x M)$ consist of forms $\om\in\Om^*(D\x M)$  which are uniquely determined  by their values $\om(\tilde X_1,\dots,\tilde X_p,Y_1,\dots,Y_q)$, where 
$X_1,\dots, X_p\in\on{Vect}_cM$ and $Y_1,\dots,Y_q\in\on{Vect}M$. 
and the exterior derivative $d=\{d^{pq}\}$ on $\Om^*(D\x M)$ is uniquely determined by the following standard formulas:
\begin{multline}\label{D1}
(d^{pq}\om)(\tilde X_1,\dots,\tilde X_{p+1};Y_1,\dots,Y_q)\\=
\sum_{i=1}^{p+1}(-1)^{i-1}\Bigl(\tilde X_i
\om(\tilde X_1,\dots,\widehat{\tilde X_i},\dots,\tilde X_{p+1};Y_1,\dots,Y_q)\\
-\sum_{j=1}^q\om(\tilde X_1,\dots,\widehat{\tilde X_i},\dots,\tilde X_{p+1};Y_1,\dots,[X_i,Y_j],\dots,Y_q)\Bigr)\\
+\sum_{i<j}(-1)^{i+j}\om([\tilde X_i,\tilde X_j],\dots,\widehat{\tilde X_i},\dots,
\widehat{\tilde X_j},\dots,\tilde X_{p+1};Y_1,\dots,Y_q),
\end{multline}
\begin{multline}\label{D2}
(d^{pq}\om)(\tilde X_1,\dots,\tilde X_p;Y_1,\dots,Y_{q+1})\\= 
(-1)^p\Bigl(\sum_{i=1}^q(-1)^{i-1}Y_i\om(\tilde X_1,\dots,\tilde X_p;Y_1,\dots,\widehat{Y_i},\dots,Y_{q+1})\\+
\sum_{i<j}(-1)^{i+j}\om(\tilde X_1,\dots,
\tilde X_p;[Y_i,Y_j],Y_1,\dots,\widehat{Y_i},\dots,\widehat{Y_i},\dots,Y_{q+1})
\Bigr),
\end{multline}
where $\om\in\Om^{pq}(G\x M)$, $X_1,\dots,X_{p+1}\in\on{Vect}_cM$ and $Y_1,\dots,Y_{q+1}\in\on{Vect}M$.

A form $\om\in\Om^{pq}(D\x M)$ is diagonal if the function  
$\om(\tilde X_1,\dots,\tilde X_p;Y_1,\dots,Y_q)$ at $(g,x)\in D\x M$ vanishes whenever the germ  of at least one of the vector fields $X_1,\dots,X_p$ at $x$ equals $0$. Denote by 
$\Om^{pq}_\De(D\x M)$ the set of diagonal $(p,q)$-forms and put 
$\Om^m_\De(D\x M)=\oplus_{p+q=m}\Om^{pq}_\De(D\x M)$. It is easy to check that 
$\Om^*_\De(D\x M)=\{\Om^m_\De(D\x M)\}$ is a subcomplex of $\Om^*(D\x M)$. 

Consider an action of the group $D$ on $D\x M$ induced by its action on $D$ by right translations and the trivial action on $M$. This action induces a left action of $D$ on
$\Om^*(D\x M)$. 
By definition, for $\om\in \Om^{pq}_\De(D\x M)$ and $g,h\in D$, we have 
$$
(h\om)(\tilde X_1,\dots,\tilde X_p;Y_1,\dots,Y_q)(g,x)=
\om(\tilde X_1^r,\dots,\tilde X_p^r;Y_1,\dots,Y_q)(h\o g,x).
$$
It is easy to see that this action of $D$ preserves the bigrading of $\Om^*_\De(D\x M)$ and the subcomplex $\Om^*_\De(D\x M)$ is $D$-stable. Denote by $H^*_\De(D\x M)$ the cohomology of the complex $\Om^*_\De(D\x M)$ and by $\Om^*_\De(D\x M)^D$ the subcomplex of $D$-invariant forms of the complex $\Om^*_\De(D\x M)$. 

\begin{lemma}\label{iso}
There is a natural isomorphism of complexes 
$$
\Om^*_\De(D\x M)^D=C^*(\on{Vect}M;\Om^*(M)).
$$
\end{lemma}
\begin{proof}
Note that, for $\om\in \Om_\De^{pq}(D\x M)^D$ and $X,X_1,\dots,X_p,Y_1,\dots,X_q$, we have 
$$
\tilde X\om(\tilde X_1,\dots,\tilde X_p;Y_1,\dots,Y_q)=
X\om(\tilde X_1,\dots,\tilde X_p;Y_1,\dots,Y_q).
$$  
Consider an action of the Lie algebra $\widetilde{\on{Vect}M}$ on $\Om^q(M)$ as follows: for $\tilde X\in\widetilde{\on{Vect}M}$ and $\om\in\Om^q(M)$, 
put $\tilde X\om=\on{L}_X\om$. 
 
Recall the standard formula for the Lie derivative $\on{L}_X\om$ of the form $\om\in\Om^q(M)$ along a vector field $X\in\on{Vect}M$.
\begin{equation}\label{Lieder}
(\on{L}_X\om)(Y_1,\dots,Y_q)=X\om(Y_1,\dots,Y_q)-
\sum_{i=1}^q\om(Y_1,\dots,[X,Y_i],\dots,Y_q),
\end{equation}
where $Y_1,\dots,Y_q\in\on{Vect}M$. 

Evidently, the form $\om\in \Om^{pq}_\De(D\x M)^D$ can be considered as a cochain of the complex 
$C^{pq}_\De(\widetilde{\on{Vect}M},\Om^q(M))$. Then, for $\om\in C_\De^{pq}(D\x M)^D$, 
using  \eqref{D1}, \eqref{D2}, and \eqref{Lieder}, we get   
\begin{multline}\label{d1}
(d^{pq}\om)(\tilde X_1,\dots,\tilde X_{p+1};\cdots)\\=
\sum_{i=1}^{p+1}(-1)^{i-1}\on{L}_{X_i}\om(\tilde X_1,\dots,\widehat{\tilde X_i},
\dots,\tilde X_{p+1};\cdots)\\
+\sum_{i<j}(-1)^{i+j}\om([\tilde X_i,\tilde X_j],\dots,\widehat{\tilde X_i},\dots,
\widehat{\tilde X_j},\dots,\tilde X_{p+1};\cdots), 
\end{multline}
\begin{multline}\label{d2}
(d^{pq}\om)(\tilde X_1,\dots,\tilde X_p;Y_1,\dots,Y_{q+1})\\= 
(-1)^p\Bigl(\sum_{i=1}^{q+1}Y_i
\om(\tilde X_1,\dots,\tilde X_p;Y_1,\dots,\widehat{Y_i},\dots,Y_{q+1})\\+
\sum_{i<j}(-1)^{i+j}\om(\tilde X_1,\dots,\tilde X_p;[Y_i,Y_j],
Y_1,\dots,\hat{Y_i},\dots,\hat{Y_j},\dots,Y_{q+1})
\Bigr).
\end{multline}

Since $X\mapsto\tilde X$ is an isomorphism $\widetilde{\on{Vect}M}=\on{Vect}M$, comparing the differential of the complex $C^*_\De(\on{Vect}M;\Om^*(M))$ with the differential of the complex $\Om^*_\De(D\x M)^D$ defined by formulas \eqref{d1} and \eqref{d2}, we see that these complexes are naturally isomorphic. 
\end{proof}

\begin{lemma}\label{Hp} For any $m,p,q$ we have 
\begin{enumerate}
\item $H^m_{\on{cont}}(D,\Om_\De^{pq}(D\x M))=0$ for $m>0$;
\item $H^0_{\on{cont}}(D,\Om^{pq}_\De(D\x M))=C^p_\De(\on{Vect}M;\Om^q(M))$.
\end{enumerate}
\end{lemma} 
\begin{proof}
Define the standard operator $B=\{B^{pq}_m\}$, where
$$
B^{pq}_m:C^m_{\on{cont}}(D,\Om^{pq}(D\x M))\to 
C^{m-1}_{\on{cont}}(D,\Om^{pq}(D\x M)), 
$$
as follows:
for $m>0$ and $c\in C^m_{\on{cont}}(D,\Om^{pq}(D\x M))$, put  
$$
(B^{pq}_mc)(g_1,\dots,g_{m-1})(\cdot)(g,x)=c(g,g_1,\dots,g_{m-1})(\cdot)(e,x),
$$
where $g,g_1,\dots,g_{m-1}\in D$ and $x\in M$. For $p=0$, put $B^{0q}=0$.

It is easy to check that for $m>0$ we have  
$$
d^{m-1}\o B^{pq}_m+B_{m+1}^{pq}\o d^m=\on{id}.
$$ 
For $c\in C^0_{\on{cont}}(D,\Om^{pq}(D\x M))$, we have 
$$
(B_1^{pq}\o d^0)(c)(\cdot)(g,x)=c(\cdot)(g,x)-c(\cdot)(e,x).
$$
Thus, $B=(B^{pq}_m)$ is the homotopy operator of the identity map of the complex  $C^*_{\on{cont}}(D,\Om^{pq}(D\x M))$ and the map which is trivial on the cochains of positive dimension and is equal 
$c(\cdot)(g,x)\mapsto c(\cdot)(e,x)$ on $C^0_{\on{cont}}(D,\Om^{pq}(D\x M))$. Later we identify 
$c(\cdot)(e,x)$ with the cochain $\bar c\in C^0_{\on{cont}}(D,(\Om^{pq}(D\x M))^D$ given by 
$\bar c(\cdot)(g,x)=c(\cdot)(e,x)$.

Since 
$B(C^*_{\on{cont}}(D,\Om^{pq}_\De)(D\x M))\subset C^*_{\on{cont}}(D,\Om^{pq}_\De(D\x M))$, the corresponding statements are true for $C^*(D,\Om^*_\De(D\x M))$. Note that one can consider 
the map $c(\cdot)(g,x)\mapsto c(\cdot)(e,x)$ on 
$C^0_{\on{cont}}(D,\Om^{pq}_\De(D\x M))^D))$ as the map 
$$C^0_{\on{cont}}(D,\Om^{pq}_\De(D\x M)^D)\to 
\Om^{pq}_\De(D\x M)^D=C^p_\De(\on{Vect}M;\Om^q(M)).
$$ 
These statements imply the claims of the lemma.
\end{proof}
By lemma \eqref{iso}, we have the following composition of homomorphisms of complexes 
\begin{multline*}
C^*_\De(\on{Vect}M;\Om^*(M))=\Om^*_\De(D\x M)^D\subset\\
 C^0_{\on{cont}}(D;\Om^*_\De(D\x M))\subset C^*_{\on{cont}}(D;\Om^*_\De(D\x M)).
\end{multline*}

\begin{corollary}\label{inclusion}
The inclusion $C^*_\De(\on{Vect}M;\Om^*(M))\subset C^*_{\on{cont}}(D;\Om^*_\De(D\x M))$
induces an isomorphism $H^*_\De(\on{Vect}M;\Om^*(M))=H^*(D;\Om^*_\De(D\x M))$.
\end{corollary}
\begin{proof}
Consider the second filtration of the double complex $C^*(D;\Om^*_\De(D\x M))$ and the corresponding 
spectral sequence $E_{2,r}$. By lemma \ref{Hp} we have $E_{2,1}^{p\,q}=0$ for $p>0$ and 
$E_{2,1}^{0\,q}=C^q_\De(\on{Vect}M;\Om^*(M))$. It is evident that the differential 
$d_1:E_{2,1}^{0\,q}\to E_{2,1}^{0\,q+1}$ coincides with the differential of the complex 
$C^*_\De(\on{Vect}M;\Om^*(M))$. Then $E_{2,2}^{p\,q}=E_{2,\infty}^{p\,q}=0$ for $p>0$ and 
$E_{2,2}^{0\,q}=E_{2,\infty}^{0\,q}=H^q_\De(\on{Vect}M;\Om^*(M))$. 

This implies the statement of the corollary.  
\end{proof}

\subsection{ A filtration of $C^*_\De(\on{Vect}M;\Om^*(M))$  and the corresponding spectral 
sequence}\label{spectral sequence}

Denote by $\Om_n^*$ the $DG$-algebra of formal differential forms in $n$ variables, i.e. the $DG$-algebra of $\infty$-jets at $0\in\R^n$ of differential forms on $\R^n$. It is clear that $\Om_n^*$ is a $W_n$-module with respect to
the action of $W_n$ by the formal Lie derivatives $\on{L}_\xi$, where $\xi\in W_n$.  Consider the complex 
$C^*(W_n,\Om^*_n)$ of standard cochains of $W_n$ with values in $\Om^*_n$ and endow it with the second differential induced by the formal exterior derivative on $\Om^*_n$. We denote by $C^*(W_n;\Om^*_n)$ the corresponding $DG$-algebra with respect to the total differential $D$ and the total grading. For $f^i\in\R[[\R^n]]$ and $\xi=\sum_{i=1}^nf^i\frac{\p}{\p x^i}\in W_n$, put 
$$
f^i_{j_1\dots j_r}(\xi)=\frac{\p^r\xi^i}{\p x^{j_1}\dots\p x^{j_r}},
$$
where $x^i$ for $i=1,\dots,n$ are the standard coordinates in $\R^n$. It is clear that $f^i_{j_1\dots j_r}\in C^1(W_n;\Om^0_n)$. Moreover, $f^i_{j_1\dots j_r}$ and $dx^i$ for $r=0,1,\dots$ and $i,j_1\dots j_r=1,\dots,n$ are generators of the $DG$-algebra $C^*(W_n;\Om_n)$.   
Since $D$ is an antiderivation of degree 1 of $DG$-algebra $C^*(W_n;\Om^*_n)$, it is uniquely determined by the following conditions: 
$$
Df^i_{j_1\dots j_r}=
\sum_{1\le k\le r}\sum_{s_1<\dots s_k}\sum_{l=1}^nf^i_{lj_1\dots\widehat{j_{s_1}}\dots
\widehat{j_{s_k}}\dots j_r}\wedge f^l_{j_{s_1}\dots j_{s_k}}
-\sum_{l=1}^nf^i_{lj_1\dots j_r}\wedge dx^l 
$$
and $D(dx^i)=f^i_j\wedge dx^i$. 

Consider the $DG$-algebra $C^*(W_n)\ox\La^*((\R^n)')$, where $(\R^n)'$ is a dual vector space for $\R^n$, with the differential which equals the differential of the complex $C^*(W_n)$ on the first factor and trivial on the second factor.

Consider the generators $c^i_{j_1,\dots,j_r}$ of the $DG$-algebra $C^*(W_n)$ as cochains of the complex $C^*(W_n)\ox\La^*((\R^n)')$. Then $c^i_{j_1,\dots,j_r}$ and 
$dx^i$ are the generators of the $DG$-algebra $C^*(W_n)\ox\La^*((\R^n)')$. 

Consider a morphism $\mu:C^*(W_n;\Om^*_n)\to C^*(W_n)\ox\La^*((\R^n)')$ of graded algebras defined by the following conditions:
$$
\mu(f^i)=dx^i+c^i,\quad \mu(f^i_{j_1,\dots,j_r})=c^i_{j_1\dots j_r}\quad\text{and}\quad \mu(dx^i)=-c^i, 
$$ 
where $r=1,2,\dots$ and $i,j_1,\dots,j_r=1,\dots,n$. 
\begin{lemma}\label{formal}
The morphism $\mu$ is an isomorphism of $DG$-algebras.
\end{lemma}
\begin{proof}
It is clear that $\mu$ is an isomorphism of graded algebras. 
The formulas above for the differential $D$ of the complex $C^*(W_n;\Om^*)$  and \eqref{dc} imply that $\mu$ is an isomorphism of $DG$-algebras. \end{proof} 

Consider again the $DG$-algebra $C^*_\De(\on{Vect}M;\Om^*(M))$. First assume that $M=U$ is a connected open subset of $\R^n$. For 
$X=\sum_{i=1}^nX^i\frac{\p}{\p x^i}\in\on{Vect}U$, where $X^i$ is a smooth function on $U$, put 
$$
\vh^i_{j_1\dots j_r}=\frac{\p^rX^i}{\p x^{j_1}\dots\p x^{j_r}}.
$$
By definition, we have $\vh^i_{j_1\dots j_r}\in C^1_\De(\on{Vect}U;\Om^0(U))$. 
By definition, the $DG$-algebra $C^*_\De(\on{Vect}U;\Om^*(U))$ is generated by the ring of smooth functions on $U$, the chains $\vh^i_{j_1\dots j_r}$ for $i,j_1,\dots,j_r=1,\dots,n$ and $r\ge 0$, and the forms $dx^i+\vh^i$ for $i=1,\dots,n$.
Then the differential $d$ of the complex $C^*_\De(\on{Vect}U;\Om^*(U))$ is uniquely defined by the following conditions: 
\begin{enumerate}
\item For a smooth function $f$ on $U$ we have 
$df=\sum_{i=1}^n\frac{\p f}{\p x_i}(dx^i+\vh^i)$;
\item For $r\ge 0$ we have 
\begin{multline*}
d\vh^i_{j_1\dots j_r}=\sum_{0\le k\le r}\sum_{s_1<\dots s_k}\sum_{l=1}^n\vh^i_{lj_1\dots\widehat{j_{s_1}}\dots
\widehat{j_{s_k}}\dots j_r}\wedge\vh^l_{j_{s_1}\dots j_{s_k}}\\
-\sum_{l=1}^n\vh^i_{lj_1\dots j_r}\wedge (dx^l+\vh^l);
\end{multline*}
\item $d(dx^i+\vh^i)=0$. 
\end{enumerate}

It is clear that each cochain $c\in C^p_\De(\on{Vect}U;\Om^q(U))$ is uniquely represented as a sum of the cochains of the following types:
\begin{equation}\label{order}
f\vh^{i_1}_{j_1\dots j_{r_1}}\wedge\dots\wedge\vh^{i_r}_{j_{r_{p-1}+1}\dots j_{r_p}}
\wedge dx^{k_1}\wedge\dots\wedge dx^{k_q},
\end{equation}
where $f\in\Om^0(U)$ and $0\le r_1\le\dots\le r_p$. 

If the germ at $x\in U$ of the cochain \eqref{order} is nonzero the number $r_p$ is called the order of the cochain \eqref{order} at $x$.  
If $c\in C^p_\De(\on{Vect}M;\Om^q(M))$ is a sum of cochains of types \eqref{order},  the maximal order of the summands at $x$ is called the order of $c$ at $x$ and is denoted by $\on{ord}_x(c)$. By definition, we have 
$\on{ord}_x(c_1+c_2)\le \max\{\on{ord}_x(c_1),\on{ord}_x(c_2)\}$. It is easy to see that $\on{ord}_x(c)$ is independent of the choice of coordinates in $U$. 

Let $M$ be an arbitrary manifold and let $c\in C^p_\De(\on{Vect}M;\Om^q(M))$ and $c\ne 0$.  Then $\on{ord}(c)$ is the maximum of orders of $c$ for any $x\in M$. If $c=0$, we put 
$\on{ord}(c)=\-\infty$.
  
Define a filtration $F_s$ of the complex $C^*_\De(\on{Vect}M;\Om^*(M))$ as follows: 
let $c\in C^p_\De(\on{Vect}M;\Om^q(M))$, then $c\in F_s$ if $\on{ord}(c)\le p+q-s$. It is clear that $F_s$ is a subalgebra of the 
differential algebra $C^*_\De(\on{Vect}M;\Om^*(M))$, $F_{s_1}\wedge F_{s_2}\subset F_{s_1+s_2}$, and $F_s\subset F_{s-1}$. 

Consider the spectral sequence $E_r$ induced by this filtration. First consider the term $E_0$ for $M=U\subset\R^n$. By definition, we can identify $E_0^{pq}$ with the vector space generated by  cochains \eqref{order} of order $s$ and the differential $d_0:E_0\to E_0$ is uniquely determined by the following conditions:
\begin{enumerate}
\item For a smooth function $f$ on $M$, we have $d_0f=0$;
\item For $r\ge 0$ we have 
\begin{multline*}d_0\vh^i_{j_1\dots j_r}=
\sum_{0\le k\le r}\sum_{s_1<\dots s_k}\sum_{l=1}^n\vh^i_{lj_1\dots\widehat{j_{s_1}}\dots
\widehat{j_{s_k}}\dots j_r}\wedge\vh^l_{j_{s_1}\dots j_{s_k}}\\
-\sum_{l=1}^n\vh^i_{lj_1\dots j_r}\wedge(dx^l+\vh^l);
\end{multline*}
\item $d_0(dx^i+\vh^i)=0$.  
\end{enumerate}
Then, comparing the differentials of the generators $c^i$,  $c^i_{j_1,\dots,j_r}$, 
and $dx^i$ of the complex $C^*(W_n)\ox\La^*((\R^n)')$ with the differentials of the generators $\vh^i_{j_1\dots j_r}$ and $dx^i+\vh^i$ of the complex $(E_0,d_0)$, we see that the complex $(E_0,d_0)$ is isomorphic to the tensor product of the complex complex 
$C^*(W_n;\Om^*_n)$ and the complex $\Om^0(U)$ with the trivial differential, i.e. we get
$$
E_1=\Om^0(U)\ox H^*(W_n)\ox\La^*((\R^n)'),
$$ 
where $dx^i+\vh^i$ are generators of the exterior algebra $\La^*((\R^n)')$. 

Let $M$ be an arbitrary manifold. Using the partition of unity on $M$ one can prove that 
$E_1=\Om^*(M)\otimes H^*(W_n)$. It is clear that the differential $d_1:E_1\to E_1$ is trivial on $H^*(W_n)$ and is equal to the exterior derivative on 
$\Om^*(M)$. Then we have $E_2=H^*(M)\otimes H^*(W_n)$.
 
\begin{remark}\label{odd} The differentials $d_r$ for $r\ge 2$ of this spectral sequence are calculated in \cite{L1} and \cite{L2}. Since the odd Pontrjagin classes of $M$ are trivial, from this calculation it follows that  for $m>2n$ all elements of $H^m(W_n,\on{O}(n))\subset H^m(W_n)\subset E_2$ live in $E_\infty$, i.e. $H^m(W_n,\on{O}(n))\subset E_\infty$.
\end{remark}
 
\subsection{The proof of main results} 

Consider a Riemannian metric on $M$ and the corresponding bundle $O(M)$ of orthogonal frames of $M$. Let $\si:O(M)\to S(M)$ be a map defined in \ref{S(M)}. 

Consider the map $f_\si:D\x O(M)\to S(M)$ given by $(g,r)\mapsto g(\si(r))$, 
where $r\in O(M)$ and $g\in D$. Consider the actions of $D$ on $D\x O(M)$ induced by the action of $D$ on $D$ by right translations and the trivial action of $D$ on $M$. Then  
the map $f_\si$ is $D$-equivariant. Moreover, the map $\si$ is $\on{O}(n)$-equivariant with respect to the natural actions of the group $\on{O}(n)$ on $O(M)$ and $S(M)$. Then the map $\si$ induces a smooth map $\tilde\si:M=O(M)/\on{O}(n)\to S(M)/\on{O}(n)$. Denote by $f_{\tilde\si}$ the map $D\x M\to S(M)/\on{O}(n)$ induced by the map $f_\si$. Consider the map $f^*_{\tilde\si}:\Om^*(S(M)/\on{O}(n))\to\Om^*(D\x M)$ 
induced by $f_{\tilde\si}$. Since the actions of $D$ and $\on{O}(n)$ on $D\x O(M)$ commute, the map $f^*_{\tilde\si}$ is $D$-equivariant homomorphism of complexes. It is easy to see 
that $f^*_{\tilde\si}(\Om^*(S(M)/\on{O}(n)))\subset \Om^*_\De(D\x M)$.  
Then we get a homomorphism of double complexes 
$$
C^*_{\on{cont}}(D;\Om^*(S(M)/\on{O}(n)))\to C^*_{\on{cont}}(D;\Om^*_\De(D\x M))
$$
induced by $f^*_{\tilde\si}$ and, by corollary \ref{inclusion}, the corresponding cohomology homomorphism 
\begin{equation}\label{hom}
H^*_{\on{cont}}(D;\Om^*(S(M)/\on{O}(n)))\to H^*_\De(\on{Vect}M;\Om^*(M)).
\end{equation}
Moreover, we have the homomorphism 
$$
\Om^*(S(M)/\on{O}(n))^D \to C^*_{\on{cont}}(D;\Om^*(S(M)/\on{O}(n))
$$
and, hence, the homomorphism 
$$
F_\si:C^*(W_n,\on{O}(n))=\Om^*(S(M)/\on{O}(n))^D\to C^*_\De(\on{Vect}M;\Om^*(M)).
$$

Now we prove the main theorem. 
\begin{theorem}\label{main}
For $m>2n$, the map $H^m(W_n,\on{O}(n))\to H^*(\on{Vect}M;\Om^*(M))$ induced by $F_\si\o\al$ is a monomorphism. 
\end{theorem}
\begin{proof}
Consider the homomorphism of complexes 
$$
f^*_{\tilde\si}:\Om^*(S(M)/\on{O}(n))^D\to\Om^*_\De(D\x M)^D.
$$ 
It will be convenient to us to treat 
$\Om^*_\De(D\x M)^D$ as a subcomplex of the de Rham complex $\Om^*(D\x O(M))$. Moreover, we will consider the forms of $\Om^*_\De(D\x O(M))$ as skew-symmetric multilinear 
forms on the Lie algebra of vector fields on $D\x O(M)$. Let $\bar X$ be the horizontal lift of a vector field $X\in\on{Vect}M$ with respect to the Levi-Civita connection. Put $\overline{\on{Vect}}M=\{\bar X;X\in\on{Vect}M\}$. Since we are interested only in forms from $\Om^*_\De(D\x M)^D$, it suffices to us to consider these forms only as multilinear functions on the vector space $\on{Vect}_cM\oplus\overline{\on{Vect}}M$. 
 
Let $(g,r)\in D\x O(M)$. Consider the linear map 
$$
(f_{\si})_*:T_{(g,r)}(D\x O(M))\to T_{\si(g,r)}S(M).
$$
For $X\in\on{Vect}M$, denote by $\tilde X$ a vector field on $S(M)$ induced by $X$. Then we have 
\begin{equation}\label{tang1}
(f_\si)_*(X^r)(g,r))=g_*\tilde X(\si(r)),
\end{equation}
where $X\in\on{Vect}_cM$, $g\in D$, and $r\in O(M)$.  For the Gelfand-Kazhdan form $\om$, by ref{G-K} we get
\begin{equation}\label{tang2}
f_\si^*\om(X^r)(g,r)=(g^*\om)(\tilde X(\si(r))=\om(\tilde X)(\si(r)). 
\end{equation}

Recall that $\om^i_{j_1\dots j_r}=\al(c^i_{j_1\dots j_r})$. 
Denote by $\nabla_i$ the operator of covariant derivative with respect to the Levi-Civita connection. It is known that, for $X\in\on{Vect}M$, $r\in \si(M)$, and the geodesic coordinates $x^i$ corresponding to $r$, we have 
$\frac{\p X^i}{\p x^j}(0)=(\nabla_jX^i)(0)$. Then, by \eqref{om(X)}, on $\si(M)$ in any coordinates on $M$ we have
\begin{equation}\label{Ga}
\om^i_j(\tilde X)=-\nabla_jX^i=
-\left(\frac{\p X^i}{\p x^j}+\sum_{k+1}^n\Gamma^i_{jk}X^k\right),
\end{equation}
where $\Gamma^i_{jk}$ are the Christoffel symbols for the Levi-Civita connection.

Define two 1-forms $\nabla^i$ and $\nabla^i_j$ on $D\x O(M)$ by the equations $\nabla^i(X^r)(g,r)=X^i(r)$, $\nabla^i(Y)(g,r)=0$, $\nabla^i_j(X^r)(g,r)=\nabla_jX^i(r)$, and $\nabla^i_j(Y)(g,r)=0$ where $X\in\on{Vect}_cM$, $Y\in\on{Vect}M$, $X^i(r)$ and $\nabla_jX^i(r)$ are the components of $X$ and its covariant derivative with respect to the frame $r\in O(M)$.
Then, by the definition of covariant derivative, \eqref{con}, \eqref{om(X)} and \eqref{tang2}, we get 
\begin{gather}
(f_\si^*\om^i)(X^r+\bar Y)(g,r))=
-\nabla^i(X^r)(g,\si(r))-\theta^i(\bar Y)(\si(r)), \label{i} \\ 
(f_\si^*\om^i_j)(X^r+\bar Y)(g,r)=-\nabla^i_j(X^r)(\si(r)),  \label{ij}
\end{gather}
where $X\in\on{Vect}_cM$, $Y\in\on{Vect}M$, and we identify the form $\theta^i$ on $O(M)$ with the corresponding form on $D\x O(M)$.  Since 
$\theta^i(\bar Y)(\si(r))=Y^i(\si(r))$, we get 

\begin{gather}
(f_\si^*\om^i)(\hat X+\bar Y)=-(\theta^i+\nabla^i)(\hat X+\bar Y),\label{I}\\ 
(f_\si^*\om^i_j)(\hat X+\bar Y)=-\nabla^i_j(\hat X+\bar Y),\label{IJ}
\end{gather}
where $\hat X=-X^r+\bar X$. 

Consider the cochains $\psi^i_j\in C^2(W_n,\on{GL}_n(\R))$ and $\Psi=(\psi^i_j)$ introduced in  \ref{S(M)} and put $\bar\Psi=\al(\Psi)$. 
Define a 1-form  $\nabla^i_{jk}$ on $D\x O(M)$ by the equations 
$\nabla^i_{jk}(X^r)(g,r)=\nabla_j\nabla_k X^i(\si(r))$ and $\nabla^i_{jk}(Y)(g,r)=0$, where $X\in\on{Vect}_cM$ and $Y\in\on{Vect}M$. By the definition of covariant derivative and the properties of the curvature tensor of a Riemannian manifold, one could check that we have 
\begin{multline}\label{ijk}
f_\si^*\bar\Psi(\hat X_1+\bar Y_1,\hat X_2+\bar Y_2)\\=-\left(\sum_k\nabla^i_{jk}\wedge(\theta^k+\nabla^k)-
\sum_{k,l}R^i_{jkl}\nabla^k\wedge\nabla^l\right)(\hat X_1+\bar Y_1,\hat X_2+\bar Y_2),
\end{multline} 
where $R^i_{jkl}$ are the components of the curvature tensor of the Levi-Civita connection.

Put $C_{p_1\dots p_l,r_1\dots r_k}=\al(c_{p_1\dots p_l,r_1\dots r_k})$, where 
$c_{p_1\dots p_l,r_1\dots r_k}$ is the basic cocycle of the complex $C^*(W_n,\on{O}(n))$ defined by \eqref{basis1}. 
Consider the cocycles $f^*_{\tilde\si}(C_{p_1\dots p_l,r_1\dots r_k})$ of the complex $\Om^*_\De(D\x M)^D=C^*_{\on{cont}}(\on{Vect}M;\Om^*(M))$. Since $\si^*\theta ^i=dx^i$, by \eqref{Ga}, \eqref{i}, \eqref{ij}, and \eqref{ijk}, the leading term (with respect to the order of cochains) of $f^*_{\tilde\si}(C_{p_1\dots p_l,r_1\dots r_k})$ equals the cocycle 
$C_{p_1\dots p_l,r_1\dots r_k}\in C^*(W_n,\on{O}(n))\subset E_0$, where 
$E_0=\Om^0(M)\otimes\La^*((\R^n)')\otimes H^*(W_n,\on{O}(n))$ is the zero term of the spectral sequence studied in \ref{spectral sequence}. It is easy to see that $f^*_{\tilde\si}(C_{p_1\dots p_l,r_1\dots r_k})$ is a basic cocycle of the second term $E_2$ of the spectral sequence. By the remark \ref{odd}, for $m>2n$ we have  $H^m(W_n,\on{O}(n))\subset E_\infty$. Hence, for $m>2n$, the map $F_\si\o\al$ induces a monomorphism of $H^m(W_n,\on{O}(n))$ into $H^*(\on{Vect}M;\Om^*(M))$. 
\end{proof} 
Theorems \ref{main} and \ref{L-G}, corollaries \ref{cocycle1} and \ref{cocycle2} imply immediately the following corollaries.
\begin{corollary}\label{main1}
Assume that $H^p(M)=0$ for $p>0$. For each $m>2n$ the map $F_\si\o\al$ induces a monomorphism of $H^m(W_n,\on{O}(n))$ into $H^m(D)$. In particular,  $c(C_{p_1\dots p_l,r_1\dots r_k})$ is a nontrivial $m$-cocycle of the complex $C^*(D)$ for $m=2(p_1+\dots+p_l+r_1+\dots+r_k)-l$. 
\end{corollary}
\begin{corollary}\label{main2}
Assume that $M$ is a closed oriented manifold. For each $m>2n$ the map $F_\si\o\al$ induces a monomorphism of $H^m(W_n,\on{O}(n))$ into $H^{m-n}(D_+)$. In particular,  $c(C_{p_1\dots p_l,r_1\dots r_k})$ is a nontrivial $(m-n)$-cocycle of the complex $C^*(D_+)$  for $m=2(p_1+\dots+p_l+r_1+\dots+r_k)-l$. 
\end{corollary}

There is a problem to find explicit expressions for the cocycles 
$c(C_{p_1\dots p_l,r_1\dots r_k})$ in the cases of corollaries \ref{main1} and \ref{main2}. In principal, it is possible under the conditions of corollary \ref{main1} whenever the manifold $M$ is contractible and the homotopy for this contraction is given and under the conditions of corollary \ref{main2}. For this one need to use formulas \eqref{Lap} and the procedure for the constructing of the group $D$ in theorem \ref{group cocycle}. It is clear that thus we will get an expression for each cocycle 
$c(C_{p_1\dots p_l,r_1\dots r_k})$ via integrals and a Riemannian metric on $M$. 
\begin{example}
Let $M=\R$ and let $x$ be the standard coordinate on $\R$. Then $s\in S(\R)$ is an 
$\infty$-jet $j_0^\infty k$, where $k(t)$ is a regular at $0$ map $\R\to\R$. We take for the coordinates on $S(M)$ the derivatives $x_i=k^{(i)}(0)$ for $i\ge 0$. We 
put $y=x_0$, $y^1=\log|x^1|$, and $y^2=\frac{x^2}{x_0^2}$. 
Consider the cocycle of the Godbillon-Vey class $c_{1,1}$. It is easy to check that we have $$
C_{1,1}=\al(c_{1,1})=dy\wedge dy^1\wedge dy^2.
$$
Applying the procedure of theorem \ref{group cocycle} one could get 
\begin{equation}\label{G-V}
c(C_{1,1})=\int_x^{f(x)}\log|h'(g(t)|d\log|g'(t)|,
\end{equation}
where $f,g,h\in\on{Diff}\R$. It is easy to see that the cohomology class of the cocycle 
$c(C_{1,1})$ given by \eqref{G-V} in the complex $C^*_{\on{cont}}(\on{Diff}\R,\R)$ is independent of the choice of $x\in\R$.   
\end{example} 
\begin{example}
Let $M$ be a closed oriented Riemannian manifold, $\theta$ the form of the Levi-Civita connection, and $\upsilon$ the volume form on $M$. For any $g\in D_+$,  put $\xi(g)=g^*\theta-\theta$ and define the function $\mu$ by the condition $g^*\upsilon=\mu\upsilon$. 
Consider the cocycles $c_{1,s_1\dots s_k}$, where $s_1+\dots +s_k=n$. 
It is clear that the cohomology clases of the cocycles $c_{1,s_1\dots s_k}$ for all $s_1,\dots, s_k$ define a part of the cohomology $H^{2n+1}(W_n,\on{O}(n))$.  
The Bott cocycle $c(C_{1,s_1\dots s_k})$ (\cite{B}) is defined by the following  formula 
$$
c(C_{1,s_1\dots s_k})
=\int_M\on{Alt}_n\left(\log(\mu)\on{tr}(\xi_1\dots\xi_{s_1})\wedge \dots\wedge 
\on{tr}(\xi_{s_1+\dots s_{k-1}+1}\dots \xi_n)\right),
$$
where $g_1\dots,g_n\in D_+$, $\bar g_i=g_i\o\dots\o g_1$, $\xi_i=\xi(\bar g_i)$, and $\on{Alt}_n$ is the alternation operator in $1,\dots,n$.  By (\cite{L2}), the cocycle $c(C_{1,s_1\dots s_k})$, given by this formula, is obtained from the cocycle $C_{1,s_1\dots s_k}$ by the procedure of theorem \ref{group cocycle}. 
Hence, the cohomology classes of Bott's cocycles $c(C_{1,s_1\dots s_k})$ are linearly independent. In particular, all of them are nontrivial. 
\end{example}
Unfortunately, for the more complicated cocycles $C_{p_1\dots p_l,r_1\dots r_k}$,  as fine formulas as for Bott's cocycles are not known.

\end{document}